*Research Article*

# Expansions of the Solutions of the General Heun Equation Governed by Two-Term Recurrence Relations for Coefficients

T. A. Ishkhanyan,[1,2,3] T. A. Shahverdyan,[3] and A. M. Ishkhanyan [1,4]

[1]*Russian-Armenian University, H. Emin 123, 0051 Yerevan, Armenia*
[2]*Moscow Institute of Physics and Technology, Dolgoprudny, Moscow Region 141700, Russia*
[3]*Institute for Physical Research, NAS of Armenia, Ashtarak 0203, Armenia*
[4]*Institute of Physics and Technology, National Research Tomsk Polytechnic University, Tomsk 634050, Russia*

Correspondence should be addressed to A. M. Ishkhanyan; aishkhanyan@gmail.com

We examine the expansions of the solutions of the general Heun equation in terms of the Gauss hypergeometric functions. We present several expansions using functions, the forms of which differ from those applied before. In general, the coefficients of the expansions obey three-term recurrence relations. However, there exist certain choices of the parameters for which the recurrence relations become two-term. The coefficients of the expansions are then explicitly expressed in terms of the gamma functions. Discussing the termination of the presented series, we show that the finite-sum solutions of the general Heun equation in terms of generally irreducible hypergeometric functions have a representation through a single generalized hypergeometric function. Consequently, the power-series expansion of the Heun function for any such case is governed by a two-term recurrence relation.

## 1. Introduction

The general Heun equation [1–3], which is the most general second-order linear ordinary differential equation having four regular singular points, is currently widely encountered in physics and mathematics research (see, e.g., [1–14] and references therein). However, this equation is much less studied than its immediate predecessor, the Gauss hypergeometric equation, which is the most general equation having three regular singular points. A reason for the slow progress in the development of the theory is that the solutions of the Heun equation (as well as its four confluent reductions) in general are not expressed in terms of definite or contour integrals involving simpler functions [2, 3]. Furthermore, the convergence regions of power-series expansions near different singularities are rather restricted and several complications arise in studying the relevant connection problems [2, 3, 15]. Another general problematic point is that the power-series solutions of the Heun equation are governed by three-term recurrence relations between successive coefficients of expansions [1–3], instead of two-term ones appearing in the hypergeometric case [16–18]. As a result, in general the coefficients are not determined explicitly.

In the present paper, we show that there exist some particular choices of the involved parameters for which the recurrence relations governing the power-series expansions become two-term. In these cases the solution of the general Heun equation can be written either as a linear combination of a finite number of the Gauss hypergeometric functions or in terms of a single generalized hypergeometric function. This is a main result of the present paper.

Another major result we report here is that in the case of the expansions of the solutions of the Heun equation in terms of hypergeometric functions there also exist particular choices of the involved parameters for which the governing three-term recurrence relations for expansion coefficients become two-term. In these cases the coefficients are explicitly written in terms of the gamma functions.

Expansions of the solutions of the Heun equation in terms of the Gauss hypergeometric functions $_2F_1$, initiated by Svartholm [19], suggest a notable extension of the series technique. This is a useful approach applicable to many differential equations including those of more general type whose nature, outside a certain region of the extended complex plane containing only two regular singular points, is not necessary to be specified exactly. Expansions involving functions other than powers have been applied to the general and confluent Heun equations by many authors. The ordinary hypergeometric [19–25], confluent hypergeometric [26–29], Coulomb wave functions [30, 31], Bessel and Hankel functions [32], incomplete Beta and Gamma functions [33–35], Hermite functions [36, 37], Goursat and Appell generalized hypergeometric functions of two variables of the first kind [38, 39], and other known special functions have been used as expansion functions. A useful property suggested by these expansions is the possibility of deriving finite-sum solutions by means of termination of the series.

As far as the expansions of the general Heun equation in terms of the hypergeometric functions are concerned, in the early papers by Svartholm [19], Erdélyi [20, 21], and Schmidt [22], the intuitive intention was to apply hypergeometric functions with parameters so chosen as to match the Heun equation as closely as possible. For this reason they used functions of the form $_2F_1(\lambda + n, \mu - n; \gamma; z)$, which have matching behavior in two singular points, $z = 0$ and 1. These functions have the following Riemann $P$-symbol representation:

$$P \begin{pmatrix} 0 & 1 & \infty \\ 0 & 0 & \lambda + n & z \\ 1 - \gamma & 1 - \delta & \mu - n \end{pmatrix}. \quad (1)$$

Here $\lambda$ and $\mu$ may adopt several values provided $1 + \lambda + \mu = \gamma + \delta$ (see [25]). It is clear that the functions have matching characteristic exponents at $z = 0$ and $z = 1$, and their behavior does not match that of the Heun function at the third singular point of the hypergeometric equation $z = \infty$.

However, it has been shown that one can also use functions that have matching behavior at only one singular point [38]. Exploring this idea, in the present paper we discuss the hypergeometric expansions of the solutions of the Heun equation in terms of functions of the form $_2F_1(\alpha, \beta; \gamma_0 \pm n; z)$ which have matching behavior (i.e., characteristic exponents) only at the singular point $z = \infty$. These functions are presented by the Riemann $P$-symbol

$$P \begin{pmatrix} 0 & 1 & \infty \\ 0 & 0 & \alpha & z \\ 1 - (\gamma_0 \pm n) & 1 - (\delta_0 \mp n) & \beta \end{pmatrix}, \quad (2)$$

where $1 + \alpha + \beta = \gamma_0 + \delta_0$, and the parameters $\gamma_0$, $\delta_0$ are chosen so that the Fuchsian condition for the general Heun equation (3), i.e., $1 + \alpha + \beta = \gamma + \delta + \varepsilon$, is fulfilled: $\gamma_0 + \delta_0 = \gamma + \delta + \varepsilon$. Note that these functions differ also from the Jacobi-polynomials used by Kalnins and Miller whose functions can be written in terms of hypergeometric functions of the form $_2F_1(\lambda + n, \mu - n; \nu + 2n; z)$ [23].

In the present paper we discuss several expansions in terms of the mentioned hypergeometric functions. In general, the coefficients of the expansions obey three-term recurrence relations similar to those known from previous developments [19–25]. However, for certain choices of the involved parameters the recurrence relations reduce to two-term ones. In these exceptional cases the coefficients of the expansions are explicitly calculated. The result is expressed in terms of the gamma functions.

Discussing the conditions for deriving finite-sum solutions by means of termination of the presented series, we show that the termination is possible if a singularity of the Heun equation is an apparent one. Furthermore, we show that any finite-sum solution of the general Heun equation derived in this way has a representation through a single generalized hypergeometric function $_pF_q$. The general conclusion is then that in any such case the power-series expansion of the Heun function is governed by a two-term recurrence relation (obviously, this is the relation obeyed by the corresponding power-series for $_pF_q$).

## 2. Hypergeometric Expansions

The general Heun equation written in its canonical form is [1]

$$\frac{d^2u}{dz^2} + \left(\frac{\gamma}{z} + \frac{\delta}{z-1} + \frac{\varepsilon}{z-a}\right)\frac{du}{dz} + \frac{\alpha\beta z - q}{z(z-1)(z-a)}u = 0, \quad (3)$$

where the parameters satisfy the Fuchsian relation $1 + \alpha + \beta = \gamma + \delta + \varepsilon$. We introduce an expansion of this equation's solution of the form

$$u = \sum_{n=0}^{\infty} c_n \cdot {}_2F_1(\alpha, \beta; \gamma_0 + n; z) \quad (4)$$

with the involved Gauss hypergeometric functions obeying the equation

$$\frac{d^2u_n}{dz^2} + \left(\frac{\gamma_0 + n}{z} + \frac{\delta + \varepsilon + \gamma - \gamma_0 - n}{z-1}\right)\frac{du_n}{dz} + \frac{\alpha\beta}{z(z-1)}u_n = 0. \quad (5)$$

Substitution of (4)-(5) into (3) gives

$$\sum_n c_n \left[\left(\frac{\gamma - \gamma_0 - n}{z} - \frac{\varepsilon + \gamma - \gamma_0 - n}{z-1} + \frac{\varepsilon}{z-a}\right)\frac{du_n}{dz} + \frac{\alpha\beta a - q}{z(z-1)(z-a)}u_n\right] = 0 \quad (6)$$

or $\sum_n c_n \left\{\left[(a-1)(\varepsilon + \gamma - \gamma_0 - n)z - a(\gamma - \gamma_0 - n)(z-1)\right]\frac{du_n}{dz} + (\alpha\beta a - q) \cdot u_n\right\} = 0. \quad (7)$

Now, using the following relations between the involved hypergeometric functions [16–18]:

$$z \frac{du_n}{dz} = \gamma_{n-1} [u_{n-1} - u_n], \quad (8)$$

$$(z-1) \frac{du_n}{dz} = -\delta_{n+1} u_n + \left(\delta_{n+1} - \frac{\alpha\beta}{\gamma_n}\right) u_{n+1}, \quad (9)$$

where $\gamma_n = \gamma_0 + n$ and $\delta_n = \delta + \varepsilon + \gamma - \gamma_n$, this equation is rewritten as

$$\sum_n c_n \left[ (a-1)(\varepsilon + \gamma - \gamma_n)(\gamma_n - 1)[u_{n-1} - u_n] \right.$$
$$+ a(\gamma - \gamma_n)\left((\delta_n - 1)u_n - \left(\delta_n - 1 - \frac{\alpha\beta}{\gamma_n}\right)u_{n+1}\right) \quad (10)$$
$$\left. + (\alpha\beta a - q) u_n \right] = 0,$$

from which we get a three-term recurrence relation for the coefficients of the expansion (4)

$$R_n c_n + Q_{n-1} c_{n-1} + P_{n-2} c_{n-2} = 0 \quad (11)$$

with

$$R_n = (a-1)(\varepsilon + \gamma - \gamma_n)(\gamma_n - 1), \quad (12)$$

$$Q_n = -R_n + a(\gamma - \gamma_n)(\alpha + \beta - \gamma_n) + (\alpha\beta a - q), \quad (13)$$

$$P_n = \frac{a}{\gamma_n}(\gamma - \gamma_n)(\alpha - \gamma_n)(\beta - \gamma_n). \quad (14)$$

From the initial conditions $c_0 = 1$ and $c_{-1} = c_{-2} = 0$ we get $(\varepsilon + \gamma - \gamma_0)(\gamma_0 - 1) = 0$. Since $\gamma_0 = 1$ is forbidden (it causes division by zero at $n = 1$ in $P_{-1}$) we obtain that the only possibility is $\gamma_0 = \varepsilon + \gamma$. Hence, the expansion finally reads

$$u = \sum_{n=0}^{\infty} c_n \cdot {}_2F_1(\alpha, \beta; \gamma + \varepsilon + n; z), \quad (15)$$

and the coefficients of the three-term recurrence relation (11) are

$$R_n = (1-a)n(\varepsilon + \gamma + n - 1), \quad (16)$$

$$Q_n = -R_n + a(1 + n - \delta)(n + \varepsilon) + (a\alpha\beta - q), \quad (17)$$

$$P_n = -\frac{a}{n+\varepsilon+\gamma}(n+\varepsilon)(n+\varepsilon+\gamma-\alpha)(n+\varepsilon+\gamma-\beta). \quad (18)$$

The expansion applies if $\alpha, \beta$, and $\gamma + \varepsilon$ are not zero or negative integers. The restrictions on $\alpha$ and $\beta$ assure that the hypergeometric functions are not polynomials of fixed degree.

The derived expansion terminates if two successive coefficients vanish. If $c_N$ is the last nonzero coefficient and $c_{N+1} = c_{N+2} = 0$ for some $N = 0, 1, 2, \ldots$, we obtain from (11) that it should be $P_N = 0$ so that the termination is possible if

$$\varepsilon = -N$$
$$\text{or } \varepsilon + \gamma - \alpha = -N \quad (19)$$
$$\text{or } \varepsilon + \gamma - \beta = -N.$$

Note that the equation $c_{N+1} = 0$ results in a polynomial equation of degree $N + 1$ for the accessory parameter $q$. This equation is convenient for rewriting the recurrence relation (11) in the following matrix form:

$$\begin{bmatrix} Q_0 & R_1 & 0 & & \\ P_0 & Q_1 & R_2 & 0 & \\ 0 & P_1 & Q_2 & R_3 & \\ & \ddots & \ddots & \ddots & \\ & & 0 & P_{N-1} & Q_N \end{bmatrix} \begin{bmatrix} c_0 \\ c_1 \\ c_2 \\ \vdots \\ c_N \end{bmatrix} = \begin{bmatrix} 0 \\ 0 \\ 0 \\ \vdots \\ 0 \end{bmatrix}. \quad (20)$$

The vanishing of the determinant of the above matrix gives the polynomial equation for $q$ defining in general $N+1$ values for which the termination occurs.

One may consider a mirror expansion

$$u = \sum_{n=0}^{\infty} c_n \cdot {}_2F_1(\alpha, \beta; \gamma_0 - n; z), \quad (21)$$

which differs from expansion (4) only by the sign of $n$ in the lower parameter of the involved hypergeometric functions. This change of the sign leads to a three-term recurrence relation (11) with the coefficients

$$R_n = \frac{a}{\gamma_0 - n}(\gamma - \gamma_0 + n)(\alpha - \gamma_0 + n)(\beta - \gamma_0 + n), \quad (22)$$

$$Q_n = -P_n + a(\gamma - \gamma_0 + n)(\alpha + \beta - \gamma_0 + n) + \alpha\beta a - q, \quad (23)$$

$$P_n = (a-1)(\varepsilon + \gamma - \gamma_0 + n)(\gamma_0 - n - 1), \quad (24)$$

where

$$\gamma_0 = \gamma \text{ or } \alpha \text{ or } \beta. \quad (25)$$

This expansion applies if $\alpha$ and $\beta$ are not zero or negative integers and $\gamma$ is not an integer.

In order for the series to terminate at some $n = N$ we put $P_N = 0$ so that this time we derive

$$\varepsilon, \varepsilon + \gamma - \alpha$$
$$\text{or } \varepsilon + \gamma - \beta = -N \quad (26)$$

for the expansions with $\gamma_0 = \gamma$ or $\alpha$ or $\beta$, respectively. Then the equation $c_{N+1} = 0$ again gives a $(N+1)$-degree polynomial equation for those values of the accessory parameter $q$ for which the termination occurs.

## 3. Finite-Sum Hypergeometric Solutions

It is readily shown that the finite-sum solutions derived from the above two types of expansions by the described termination procedure coincide, as can be expected because of apparent symmetry. For example, consider the expansion (15) in the case $\varepsilon = -N$. The involved hypergeometric functions have the form $u_n = {}_2F_1(\alpha, \beta, \gamma - N + n; z)$. Since $n = 0, 1, \ldots, N$, we see that the set of the involved hypergeometric functions is exactly the same as in the case of the second type expansion (21) with $\gamma_0 = \gamma$ : $\{{}_2F_1(\alpha, \beta, \gamma; z), {}_2F_1(\alpha, \beta, \gamma - 1; z), \ldots, {}_2F_1(\alpha, \beta, \gamma - N; z)\}$. Furthermore, examination of (20) shows that the equation for the accessory parameter $q$ and the expansion coefficients are also the same for the two expansions. The same happens to other two cases: $\varepsilon + \gamma - \alpha = -N$ and $\varepsilon + \gamma - \beta = -N$. Thus, while different in general, the expansions (15)-(18) and (21)-(25) lead to the same finite-sum closed-form solutions.

Consider the explicit forms of these solutions examining, for definiteness, the expansion (15)-(18). An immediate observation is that because of the symmetry of the Heun equation with respect to the interchange $\alpha \longleftrightarrow \beta$ the finite-sum solutions produced by the choices $\varepsilon + \gamma - \alpha = -N$ and $\varepsilon + \gamma - \beta = -N$ are of the same form. Furthermore, by applying the formula [16]

$$
{}_2F_1(\alpha, \beta; \alpha + k; z) = (1-z)^{k-\beta} {}_2F_1(k, \alpha - \beta + k; \alpha + k; z) \quad (27)
$$

to the involved hypergeometric functions ${}_2F_1(\alpha, \beta; \alpha - N + n; z)$ or ${}_2F_1(\alpha, \beta; \beta - N + n; z)$ we see that the sum is a quasi-polynomial, namely, a product of $(1-z)^{1-\delta}$ and a polynomial in $z$. Here are the first two of the solutions:

$$
\varepsilon + \gamma - \alpha = 0,
$$
$$
q = a\gamma(\delta - 1), \quad (28)
$$
$$
u = (1-z)^{1-\delta},
$$
$$
\varepsilon + \gamma - \alpha = -1,
$$
$$
q^2 + [\alpha - 1 - a(\delta - 2 + \gamma(2\delta - 3))]q - a\gamma(\delta - 2) \quad (29)
$$
$$
\cdot (\alpha - a(1+\gamma)(\delta - 1)) = 0,
$$
$$
u = (1-z)^{1-\delta}
$$
$$
\cdot \left(1 - \frac{\alpha + 1 - \delta}{\alpha - 1}z + \frac{q - a(\alpha\beta + \varepsilon - \delta\varepsilon)}{(1-a)(\alpha - 1)}(1-z)\right). \quad (30)
$$

Note that, since $1 - \delta$ is a characteristic exponent of the Heun equation, the transformation $u = (1-z)^{1-\delta}w(z)$ results in another Heun equation for $w(z)$. Hence, the derived finite-sum solutions corresponding to the choices $\varepsilon + \gamma - \alpha = -N$ and $\varepsilon + \gamma - \beta = -N$ are generated from the polynomial solutions of the equation for $w(z)$.

More interesting is the case $\varepsilon = -N$, when the finite-sum solutions involve $N+1$ hypergeometric functions irreducible, in general, to simpler functions. The case $\varepsilon = 0$ produces the trivial result $q = a\alpha\beta$, when the Heun equation is degenerated into the hypergeometric equation with the solution $u = {}_2F_1(\alpha, \beta; \gamma; z)$. The solution for the first nontrivial case $\varepsilon = -1$ reads

$$
u = {}_2F_1(\alpha, \beta; \gamma - 1; z) + \frac{q - a\alpha\beta + a(1-\delta)}{(1-a)(\gamma - 1)} \\
\cdot {}_2F_1(\alpha, \beta; \gamma; z), \quad (31)
$$

where $q$ is a root of the equation

$$
(q - a\alpha\beta + a(1-\delta))(q - a\alpha\beta + (a-1)(1-\gamma)) \\
- a(1-a)(1+\alpha-\gamma)(1+\beta-\gamma) = 0. \quad (32)
$$

Note that the second term in (31) vanishes if $q - a\alpha\beta + a(1-\delta) = 0$ so that in this degenerate case the solution involves one, not $2 = N + 1$, terms. It is seen from (32) that this situation is necessarily the case if $a = 1/2$, $\gamma + \delta = 2$ and $\alpha$ or $\beta$ equals $\gamma - 1$. We will see that the solution in this case is a member of a family of specific solutions for which the expansion is governed by two-term recurrence relations for the coefficients.

The solutions (31) and (32) have been noticed on several occasions [40–44]. It has been shown that, for $\varepsilon = -1$, when the characteristic exponents of $z = a$ are $0, 2$ so that they differ by an integer, (32) provides the condition for the singularity $z = a$ to be *apparent* (or "simple"); that is, no logarithmic terms are involved in the local Frobenius series expansion [40–42]. In fact, the Frobenius solution in this case degenerates to a Taylor series. It has further been observed that the solution (31) can be expressed in terms of the Clausen generalized hypergeometric function ${}_3F_2$ with an upper parameter exceeding a lower one by unity [40–42]:

$$
\frac{u}{u(0)} = {}_3F_2(\alpha, \beta, e + 1; \gamma, e; z), \quad (33)
$$

where the parameter $e$ is given as

$$
e = \frac{a\alpha\beta}{q - a\alpha\beta}. \quad (34)
$$

Note that using this parameter $e$ the solution of (32) is parameterized as [41]

$$
q = a\alpha\beta\frac{1+e}{e},
$$
$$
a = \frac{e(e-\gamma+1)}{(e-\alpha)(e-\beta)}. \quad (35)
$$

We will now show that a similar generalized hypergeometric representation holds also for $\varepsilon = -2$ and for all $\varepsilon \in \mathbb{Z}, \varepsilon \neq 1$.

## 4. The Case $\varepsilon \leq -2$, $\varepsilon \in \mathbb{Z}$

For $\varepsilon = -2$ the termination equation $c_{N+1} = 0$ for the accessory parameter $q$ is written as

$$((q - a\alpha\beta)^2 + (q - a\alpha\beta)(4a - 2 - (3 + \alpha + \beta)a + \gamma) \\ + 2a(a-1)\alpha\beta) \times (q - a\alpha\beta - 2(1 + \alpha + \beta)a - 2 \\ + 2\gamma) + (q - a\alpha\beta)2a(a-1)(\alpha\beta + 1 + \alpha + \beta) \\ = 0. \tag{36}$$

The solution of the Heun equation for a root of this equation is given as

$$u = {}_2F_1(\alpha, \beta; \gamma - 2; z) + B_1 \cdot {}_2F_1(\alpha, \beta; \gamma - 1; z) + B_2 \\ \cdot {}_2F_1(\alpha, \beta; \gamma; z), \tag{37}$$

with

$$B_1 = \frac{q - a\alpha\beta + a(1-\delta)}{(1-a)(\gamma-2)}, \\ B_2 = \frac{a(1+\alpha-\gamma)(1+\beta-\gamma)}{(q - a\alpha\beta + 2(1-a)(\gamma-1))(\gamma-1)} B_1. \tag{38}$$

It is now checked that this solution is presented by the hypergeometric function ${}_4F_3$ as [42]

$$\frac{u}{u(0)} = {}_4F_3(\alpha, \beta, e+1, r+1; \gamma, e, r; z), \tag{39}$$

where $u(0) = 1 + B_1 + B_2$ and the parameters $e$, $r$ solve the equations (compare with (35))

$$q = a\alpha\beta \frac{(1+e)(1+r)}{er}, \tag{40}$$

$$a = \frac{er(2er + (e+r+1)(2-\gamma))}{\alpha\beta((1+e)^2 + (1+r)^2 - 1) + er(2er - 4 - (e+r+3)(\alpha+\beta-1))}. \tag{41}$$

It is further checked that this system of equations admits a unique solution $e$, $r$ (up to the transposition $e \longleftrightarrow r$).

The presented result is derived in a simple way by substituting the ansatz (39) into the general Heun equation and expanding the result in powers of $z$. The equations resulting in cancelling the first three terms proportional to $z^0$, $z^1$, and $z^2$ are that given by (40), (41), and (36), respectively. It is then shown that these three equations are enough for the Heun equation to be satisfied identically.

A further remark is that (36) presents the condition for the singularity $z = a$ to be apparent for $\varepsilon = -2$. This is straightforwardly verified by checking the power-series solution $u = \sum_{n=0}^{\infty} c_n(z-a)^n$ with $c_0 \neq 0$ for the neighborhood of the point $z = a$. In calculating $c_3$ a division by zero will occur, unless $q$ satisfies (36), in which case the equation for $c_3$ will be identically satisfied.

It can be checked that generalized hypergeometric representations are achieved also for $\varepsilon = -3, -4, -5$ [42]. The conjecture is that for any negative integer $\varepsilon = -N$, $N = 1, 2, 3, \ldots$ there exists a generalized hypergeometric solution of the Heun equation given by the ansatz

$$u \\ = {}_{N+2}F_{1+N}(\alpha, \beta, e_1+1, \ldots, e_N+1; \gamma, e_1, \ldots, e_N; z) \tag{42}$$

provided the singularity at $z = a$ is an apparent one. For the latter condition to be the case, the accessory parameter $q$ should satisfy a $(N+1)$-degree polynomial equation which forces the above expansions (15)-(18) and (21)-(25) to terminate at $N$th term. Note that (42) applies also for $N = 0$, that is, for $\varepsilon = 0$, for which the Heun function degenerates to the Gauss hypergeometric function $u = {}_2F_1(\alpha, \beta; \gamma; z)$ provided $q = 0$.

Finally, we note that by the elementary power change $u = (z-a)^{1-\varepsilon} w$ a Heun equation with a *positive* exponent parameter $\varepsilon > 2$ is transformed into the one with a negative parameter $2 - \varepsilon < 0$. Hence, it is understood that a similar generalized hypergeometric representation of the solution of the Heun equation can also be constructed for positive integer $\varepsilon = N$, $N = 2, 3, \ldots$. Thus, the only exception is the case $\varepsilon = 1$.

It is a basic knowledge that the generalized hypergeometric function ${}_pF_q$ is given by a power-series with coefficients obeying a two-term recurrence relation. Since any finite-sum solution of the general Heun equation derived via termination of a hypergeometric series expansion has a representation through a single generalized hypergeometric function ${}_pF_q$, the general conclusion is that in each such case the power-series expansion of the Heun function is governed by a two-term recurrence relation (obviously, by the relation obeyed by the corresponding power-series for ${}_pF_q$).

## 5. Hypergeometric Expansions with Two-Term Recurrence Relations for the Coefficients

In this section we explore if the three-term recurrence relations governing the above-presented hypergeometric expansions can be reduced to two-term ones. We will see that the answer is positive. Two-term reductions are achieved for an infinite set of particular choices of the involved parameters.

First, we mention a straightforward case which actually turns to be rather simple because in this case the Heun equation is transformed into the Gauss hypergeometric equation by a variable change. This is the case if

$$a = 1/2,$$
$$\gamma + \delta = 2 \quad (43)$$
$$\text{and } q = a\alpha\beta + a(1-\delta)\varepsilon,$$

when the coefficient $Q_n$ in (11) identically vanishes so that the recurrence relation between the expansion coefficients straightforwardly becomes two-term for both expansions (15)-(18) and (21)-(25). The coefficients of the expansions are then explicitly calculated. For instance, expansion (15) is written as

$$u = \sum_{k=0}^{\infty} \frac{(\varepsilon/2)_k ((\gamma+\varepsilon-\alpha)/2)_k ((\gamma+\varepsilon-\beta)/2)_k}{k! ((\gamma+\varepsilon)/2)_k ((1+\gamma+\varepsilon)/2)_k} \quad (44)$$
$$\cdot {}_2F_1(\alpha, \beta; \gamma+\varepsilon+2k; z),$$

where $(\ldots)_k$ denotes the Pochhammer symbol. The values $u(0)$, $u'(0)$ and $u(1)$, and $u'(1)$ can then be written in terms of generalized hypergeometric series [17, 18]. For instance,

$$u(0) = {}_3F_2\left(\frac{\gamma+\varepsilon-\alpha}{2}, \frac{\gamma+\varepsilon-\beta}{2}, \frac{\varepsilon}{2}; \frac{\gamma+\varepsilon}{2}, \frac{1+\gamma+\varepsilon}{2}; 1\right). \quad (45)$$

However, as it was already mentioned above, the case (43) is a rather simple one because the transformation

$$u(z) = z^{1-\gamma}\left(1-\frac{z}{a}\right)^{1-\varepsilon} w(4z(1-z)) \quad (46)$$

reduces the Heun equation to the Gauss hypergeometric equation for the new function $w$. The solution of the general Heun equation is then explicitly written as

$$u(z) = z^{1-\gamma}\left(1-\frac{z}{a}\right)^{1-\varepsilon}$$
$$\cdot {}_2F_1\left(\frac{1-\alpha+\delta}{2}, \frac{1-\beta+\delta}{2}; \delta; 4(1-z)z\right). \quad (47)$$

Now, we will show that there exist nontrivial cases of two-term reductions of the three-term recurrence (11) with (16)-(18) or (22)-(24). These reductions are achieved by the following ansatz guessed by examination of the structure of solutions (33) and (39):

$$c_n = \left(\frac{1}{n} \frac{\prod_{k=1}^{N+2}(a_k - 1 + n)}{\prod_{k=1}^{N+1}(b_k - 1 + n)}\right) c_{n-1}, \quad (48)$$

where, having in mind the coefficients $R_n$, $Q_n$, and $P_n$ given by (16)-(18), we put

$$a_1, \ldots, a_N, a_{N+1}, a_{N+2}$$
$$= 1 + e_1, \ldots, 1 + e_N, \gamma+\varepsilon-\alpha, \gamma+\varepsilon-\beta, \quad (49)$$
$$b_1, \ldots, b_N, b_{N+1} = e_1, \ldots, e_N, \gamma+\varepsilon \quad (50)$$

with parameters $e_1, \ldots, e_N$ to be defined later. Note that this ansatz implies that $e_1, \ldots, e_N$ are not zero or negative integers.

The ratio $c_n/c_{n-1}$ is explicitly written as

$$\frac{c_n}{c_{n-1}} = \frac{(\gamma+\varepsilon-\alpha-1+n)(\gamma+\varepsilon-\alpha-1+n)}{(\gamma+\varepsilon-1+n)n} \prod_{k=1}^{N}$$
$$\cdot \frac{e_k+n}{e_k-1+n}. \quad (51)$$

With this, the recurrence relation (11) is rewritten as

$$R_n \frac{(\gamma+\varepsilon-\alpha-1+n)(\gamma+\varepsilon-\beta-1+n)}{(\gamma+\varepsilon-1+n)n} \prod_{k=1}^{N}$$
$$\cdot \frac{e_k+n}{e_k-1+n} + Q_{n-1} + P_{n-2}$$
$$\cdot \frac{(\gamma+\varepsilon-2+n)(n-1)}{(\gamma+\varepsilon-\alpha-2+n)(\gamma+\varepsilon-\beta-2+n)} \prod_{k=1}^{N} \quad (52)$$
$$\cdot \frac{e_k-2+n}{e_k-1+n} = 0.$$

Substituting $R_n$ and $P_{n-2}$ from (16) and (18) and cancelling the common denominator, this equation becomes

$$(1-a)(\gamma+\varepsilon-\alpha-1+n)(\gamma+\varepsilon-\beta-1+n)$$
$$\cdot \prod_{k=1}^{N}(e_k+n) + Q_{n-1}\prod_{k=1}^{N}(e_k-1+n) \quad (53)$$
$$- a(\varepsilon+n-2)(n-1)\prod_{k=1}^{N}(e_k-2+n) = 0.$$

This is a polynomial equation in $n$. Notably, it is of degree $N+1$, not $N+2$, because the highest-degree term proportional to $n^{N+2}$ identically vanishes. Hence, we have an equation of the form

$$\sum_{m=0}^{N+1} A_m(a, q; \alpha, \beta, \gamma, \delta, \varepsilon; e_1, \ldots, e_N) n^m = 0. \quad (54)$$

Then, equating to zero the coefficients $A_m$ warrants the satisfaction of the three-term recurrence relation (11) for all $n$. We thus have $N+2$ equations $A_m = 0$, $m = 0, 1, \ldots, N+1$, of which $N$ equations serve for determination of the parameters $e_{1,2,\ldots,N}$ and the remaining two impose restrictions on the parameters of the Heun equation.

One of these restrictions is derived by calculating the coefficient $A_{N+1}$ of the term proportional to $n^{N+1}$ which is readily shown to be $2 + N - \delta$. Hence,

$$\delta = 2 + N. \quad (55)$$

The second restriction imposed on the parameters of the Heun equation is checked to be a polynomial equation of the degree $N+1$ for the accessory parameter $q$.

With the help of the Fuchsian condition $1+\alpha+\beta = \gamma+\delta+\varepsilon$, we have

$$\gamma + \varepsilon - \alpha - 1 = \beta - \delta = \beta - 2 - N, \quad (56)$$

$$\gamma + \varepsilon - \beta - 1 = \alpha - \delta = \alpha - 2 - N, \quad (57)$$

$$\gamma + \varepsilon - 1 = \alpha + \beta - \delta = \alpha + \beta - 2 - N, \quad (58)$$

so that the two-term recurrence relation (51) can be rewritten as ($c_0 = 1$)

$$c_n = \left( \frac{(\alpha - 2 - N + n)(\beta - 2 - N + n)}{(\alpha + \beta - 2 - N + n)n} \prod_{k=1}^{N} \frac{e_k + n}{e_k - 1 + n} \right) c_{n-1}, \quad n \geq 1. \quad (59)$$

Note that it follows from this relation, since $\alpha$, $\beta$, and $e_1, \ldots, e_N$ are not zero or negative integers, that $c_n$ may vanish only if $\alpha$ is a positive integer such that $0 < \alpha < 2 + N$ or $\beta$ is a positive integer such that $0 < \beta < 2 + N$.

Resolving the recurrence (51), the coefficients of expansion (15)-(18) are explicitly written in terms of the gamma functions as

$$c_n = \frac{\Gamma(\gamma + \varepsilon)\Gamma(n + \gamma + \varepsilon - \alpha)\Gamma(n + \gamma + \varepsilon - \beta)}{n!\Gamma(\gamma + \varepsilon - \alpha)\Gamma(\gamma + \varepsilon - \beta)\Gamma(n + \gamma + \varepsilon)} \prod_{k=1}^{N} \frac{e_k + n}{e_k}, \quad n \geq 1. \quad (60)$$

Here are the explicit solutions of the recurrence relation (11) for $N = 0$ and $N = 1$.

$N = 0$:
$$\quad (61)$$
$$\delta = 2,$$

$$q = a\gamma + (\alpha - 1)(\beta - 1), \quad (62)$$

$$c_n = \frac{\Gamma(\gamma + \varepsilon)\Gamma(n + \gamma + \varepsilon - \alpha)\Gamma(n + \gamma + \varepsilon - \beta)}{n!\Gamma(\gamma + \varepsilon - \alpha)\Gamma(\gamma + \varepsilon - \beta)\Gamma(n + \gamma + \varepsilon)}. \quad (63)$$

$N = 1$:
$$\quad (64)$$
$$\delta = 3,$$

$$q^2 - q(4 + a - 3\alpha - 3\beta + 2\alpha\beta + 3a\gamma) + (\alpha - 2)$$
$$\cdot (\alpha - 1)(\beta - 2)(\beta - 1) \quad (65)$$
$$+ a(4 + 2a - 4\alpha - 4\beta + 3\alpha\beta)\gamma + 2a^2\gamma^2 = 0,$$

$$e_1 = -q + a(1 + \gamma) - 1 + (\alpha - 1)(\beta - 1), \quad (66)$$

$$c_n = \frac{\Gamma(\gamma + \varepsilon)\Gamma(n + \gamma + \varepsilon - \alpha)\Gamma(n + \gamma + \varepsilon - \beta)}{n!\Gamma(\gamma + \varepsilon - \alpha)\Gamma(\gamma + \varepsilon - \beta)\Gamma(n + \gamma + \varepsilon)}$$
$$\cdot \frac{e_1 + n}{e_1}. \quad (67)$$

These results are readily checked by direct verification of the recurrence relation (11) with coefficients (16)-(18). We conclude by noting that similar explicit solutions can be straightforwardly derived for the expansion (21)-(25) as well.

## 6. Discussion

Thus, we have presented an expansion of the solutions of the Heun equation in terms of hypergeometric functions having the form $_2F_1(\alpha, \beta; \gamma_0 + n; z)$ with $\gamma_0 = \varepsilon + \gamma$ and expansions in terms of functions $_2F_1(\alpha, \beta; \gamma_0 - n; z)$ with $\gamma_0 = \gamma, \alpha, \beta$. For any set of parameters of the Heun equation provided that $\gamma + \varepsilon, \gamma, \alpha, \beta$ are not all simultaneously integers at least one of these expansions can be applied. Obviously, the expansions are meaningless if $\alpha\beta = 0$ since then the involved hypergeometric functions are mere constants and for the solution the summation produces the trivial result $u = 0$.

The applied technique is readily extended to the four confluent Heun equations. For instance, the solutions of the single- and double-confluent Heun equations using the Kummer confluent hypergeometric functions of the forms $_1F_1(\alpha_0 + n; \gamma_0 + n; s_0 z)$, $_1F_1(\alpha_0 + n; \gamma_0; s_0 z)$, and $_1F_1(\alpha_0; \gamma_0 + n; s_0 z)$ are straightforward. By slight modification, equations of more general type, e.g., of the type discussed by Schmidt [22], can also be considered. In all these cases the termination of the series results in closed-form solutions appreciated in many applications. A representative example is the determination of the exact complete return spectrum of a quantum two-state system excited by a laser pulse of Lorentzian shape and of a double level-crossing frequency detuning [45]. A large set of recent applications of the finite-sum expansions of the biconfluent Heun equation in terms of the Hermite functions to the Schrödinger equation is listed in [37] and references therein.

Regarding the closed-form solutions produced by the presented expansions, this happens in three cases: $\varepsilon = -N$, $\varepsilon + \gamma - \alpha = -N$, $\varepsilon + \gamma - \beta = -N$, $N = 0, 1, 2, 3, \ldots$. In each case the general Heun equation admits finite-sum solutions in general at $N + 1$ choices of the accessory parameter $q$ defined by a polynomial equation of the order of $N + 1$. The last two choices for $\varepsilon$ result in quasi-polynomial solutions, while, in the first case, when $\varepsilon$ is a negative integer, the solutions involve $N + 1$ hypergeometric functions generally irreducible to simpler functions. Discussing the termination of this series, we have shown that this is possible if a singularity of the Heun equation is an apparent one. We have further shown that the corresponding finite-sum solution of the general Heun equation has a representation through a single generalized hypergeometric function. The general conclusion suggested by this result is that in any such case the power-series expansion of the Heun function is governed by the two-term recurrence relation obeyed by the power-series for the corresponding generalized hypergeometric function $_pF_q$.

There are many examples of application of finite-sum solutions of the Heun equation both in physics and mathematics [46–56], for instance, the solution of a class of free boundary problems occurring in groundwater flow in liquid mechanics and the removal of false singular points of Fuchsian ordinary differential equations in applied mathematics [43]. Another example is the derivation of the third independent exactly solvable hypergeometric potential, after the Eckart and the Pöschl-Teller potentials, which is proportional to an energy-independent parameter and has a shape that is independent of this parameter [44]. Some other recent examples can be found in references listed in [14].

Finally, we have shown that there exist infinitely many choices of the involved parameters for which the three-term recurrence relations governing the hypergeometric expansions of the solutions of the general Heun equation are reduced to two-term ones. The coefficients of the expansions are then explicitly expressed in terms of the gamma functions. We have explicitly presented two such cases.

## Data Availability

No data were used to support this study.

## Conflicts of Interest

The authors declare that they have no conflicts of interest.

## Acknowledgments

This research has been conducted within the scope of the International Associated Laboratory IRMAS (CNRS-France & SCS-Armenia). The work has been supported by the Russian-Armenian (Slavonic) University at the expense of the Ministry of Education and Science of the Russian Federation, the Armenian State Committee of Science (SCS Grant no. 18RF-139), the Armenian National Science and Education Fund (ANSEF Grant no. PS-4986), and the project "Leading Russian Research Universities" (Grant no. FTI_24_2016 of the Tomsk Polytechnic University). T. A. Ishkhanyan acknowledges the support from SPIE through a 2017 Optics and Photonics Education Scholarship and thanks the French Embassy in Armenia for a doctoral grant as well as the Agence Universitaire de la Francophonie and Armenian State Committee of Science for a Scientific Mobility grant.

## References

[1] K. Heun, "Zur Theorie der Riemann'schen Functionen zweiter Ordnung mit vier Verzweigungspunkten," *Mathematische Annalen*, vol. 33, no. 2, pp. 161–179, 1888.

[2] A. Ronveaux, Ed., *Heun's Differential Equations*, Oxford University Press, Oxford, UK, 1995.

[3] S. Y. Slavyanov and W. Lay, *Special functions*, Oxford Mathematical Monographs, Oxford University Press, Oxford, UK, 2000.

[4] E. Renzi and P. Sammarco, "The hydrodynamics of landslide tsunamis: Current analytical models and future research directions," *Landslides*, vol. 13, no. 6, pp. 1369–1377, 2016.

[5] M. Renardy, "On the eigenfunctions for Hookean and FENE dumbbell models," *Journal of Rheology*, vol. 57, no. 5, pp. 1311–1324, 2013.

[6] M. M. Afonso and D. Vincenzi, "Nonlinear elastic polymers in random flow," *Journal of Fluid Mechanics*, vol. 540, pp. 99–108, 2005.

[7] I. C. Fonseca and K. Bakke, "Quantum Effects on an Atom with a Magnetic Quadrupole Moment in a Region with a Time-Dependent Magnetic Field," *Few-Body Systems*, vol. 58, no. 1, 2017.

[8] Q. Xie, H. Zhong, M. T. Batchelor, and C. Lee, "The quantum Rabi model: solution and dynamics," *Journal of Physics A: Mathematical and General*, vol. 50, no. 11, 113001, 40 pages, 2017.

[9] C. A. Downing and M. E. Portnoi, "Massless Dirac fermions in two dimensions: Confinement in nonuniform magnetic fields," *Physical Review B: Condensed Matter and Materials Physics*, vol. 94, no. 16, 2016.

[10] P. Fiziev and D. Staicova, "Application of the confluent Heun functions for finding the quasinormal modes of nonrotating black holes," *Physical Review D: Particles, Fields, Gravitation and Cosmology*, vol. 84, no. 12, 2011.

[11] D. Batic, D. Mills-Howell, and M. Nowakowski, "Potentials of the Heun class: the triconfluent case," *Journal of Mathematical Physics*, vol. 56, no. 5, 052106, 17 pages, 2015.

[12] H. S. Vieira and V. B. Bezerra, "Confluent Heun functions and the physics of black holes: resonant frequencies, HAWking radiation and scattering of scalar waves," *Annals of Physics*, vol. 373, pp. 28–42, 2016.

[13] R. L. Hall and N. Saad, "Exact and approximate solutions of Schrödinger's equation with hyperbolic double-well potentials," *The European Physical Journal Plus*, vol. 131, no. 8, 2016.

[14] M. Hortacsu, *Proc. 13th Regional Conference on Mathematical Physics*, World Scientific, Singapore, 2013.

[15] R. Schäfke and D. Schmidt, "The connection problem for general linear ordinary differential equations at two regular singular points with applications in the theory of special functions," *SIAM Journal on Mathematical Analysis*, vol. 11, no. 5, pp. 848–862, 1980.

[16] A. Erdélyi, W. Magnus, F. Oberhettinger, and F. G. Tricomi, *Higher Transcendental Functions*, vol. 3, McGraw-Hill, New York, NY, USA, 1955.

[17] L. J. Slater, *Generalized Hypergeometric Functions*, Cambridge University Press, Cambridge, UK, 1966.

[18] W. N. Bailey, *Generalized Hypergeometric Series*, Stechert-Hafner Service Agency, 1964.

[19] N. Svartholm, "Die Lösung der Fuchs'schen Differentialgleichung zweiter Ordnung durch Hypergeometrische Polynome," *Mathematische Annalen*, vol. 116, no. 1, pp. 413–421, 1939.

[20] A. Erdélyi, "The Fuchsian equation of second order with four singularities," *Duke Mathematical Journal*, vol. 9, pp. 48–58, 1942.

[21] A. Erdélyi, "Certain expansions of solutions of the Heun equation," *Quarterly Journal of Mathematics*, vol. 15, pp. 62–69, 1944.

[22] D. Schmidt, "Die Lösung der linearen Differentialgleichung 2. Ordnung um zwei einfache Singularitäten durch Reihen nach hypergeometrischen Funktionen," *Journal fur die Reine und Angewandte Mathematik*, vol. 309, pp. 127–148, 1979.

[23] E. G. Kalnins and J. Miller, "Hypergeometric expansions of Heun polynomials," *SIAM Journal on Mathematical Analysis*, vol. 22, no. 5, pp. 1450–1459, 1991.

[24] S. Mano, H. Suzuki, and E. Takasugi, "Analytic solutions of the Teukolsky equation and their low frequency expansions," *Progress of Theoretical and Experimental Physics*, vol. 95, no. 6, pp. 1079–1096, 1996.

[25] B. D. Sleeman and V. B. Kuznetsov, "Heun Functions: Expansions in Series of Hypergeometric Functions," http://dlmf.nist.gov/31.11.

[26] T. Kurth and D. Schmidt, "On the global representation of the solutions of second-order linear differential equations having an irregular singularity of rank one in by series in terms of confluent hypergeometric functions," *SIAM Journal on Mathematical Analysis*, vol. 17, no. 5, pp. 1086–1103, 1986.

[27] S. Mano and E. Takasugi, "Analytic solutions of the Teukolsky equation and their properties," *Progress of Theoretical and Experimental Physics*, vol. 97, no. 2, pp. 213–232, 1997.

[28] L. J. El-Jaick and B. D. Figueiredo, "Solutions for confluent and double-confluent Heun equations," *Journal of Mathematical Physics*, vol. 49, no. 8, 083508, 28 pages, 2008.

[29] T. A. Ishkhanyan and A. M. Ishkhanyan, "Expansions of the solutions to the confluent Heun equation in terms of the Kummer confluent hypergeometric functions," *AIP Advances*, vol. 4, Article ID 087132, 2014.

[30] E. W. Leaver, "Solutions to a generalized spheroidal wave equation: Teukolsky's equations in general relativity, and the two-center problem in molecular quantum mechanics," *Journal of Mathematical Physics*, vol. 27, no. 5, pp. 1238–1265, 1986.

[31] L. J. El-Jaick and B. D. Figueiredo, "Confluent Heun equations: convergence of solutions in series of Coulomb wavefunctions," *Journal of Physics A: Mathematical and General*, vol. 46, no. 8, 085203, 29 pages, 2013.

[32] B. D. B. Figueiredo, "Generalized spheroidal wave equation and limiting cases," *Journal of Mathematical Physics*, vol. 48, Article ID 013503, 2007.

[33] A. Ishkhanyan, "Incomplete beta-function expansions of the solutions to the confluent Heun equation," *Journal of Physics A: Mathematical and General*, vol. 38, no. 28, pp. L491–L498, 2005.

[34] E. S. Cheb-Terrab, "Solutions for the general, confluent and biconfluent Heun equations and their connection with Abel equations," *Journal of Physics A: Mathematical and General*, vol. 37, no. 42, pp. 9923–9949, 2004.

[35] T. A. Ishkhanyan, Y. Pashayan-Leroy, M. R. Gevorgyn, C. Leroy, and A. M. Ishkhanyan, "Expansions of the solutions of the biconfluent Heun equation in terms of incomplete Beta and Gamma functions," *Journal of Contemporary Physics (Armenian Academy of Sciences)*, vol. 51, pp. 229–236, 2016.

[36] A. Hautot, "Sur des combinaisons linéaires d'un nombre fini de fonctions transcendantes comme solutions d'équations différentielles du second ordre," *Bulletin de la Societe Royale des Sciences de Liege*, vol. 40, pp. 13–23, 1971.

[37] T. A. Ishkhanyan and A. M. Ishkhanyan, "Solutions of the bi-confluent Heun equation in terms of the Hermite functions," *Annals of Physics*, vol. 383, pp. 79–91, 2017.

[38] A. Ishkhanyan and K.-A. Suominen, "New solutions of Heun's general equation," *Journal of Physics A: Mathematical and General*, vol. 36, no. 5, pp. L81–L85, 2003.

[39] C. Leroy and A. M. Ishkhanyan, "Expansions of the solutions of the confluent Heun equation in terms of the incomplete Beta and the Appell generalized hypergeometric functions," *Integral Transforms and Special Functions*, vol. 26, no. 6, pp. 451–459, 2015.

[40] J. Letessier, G. Valent, and J. Wimp, "Some differential equations satisfied by hypergeometric functions," in *Approximation and computation (West Lafayette, IN, 1993)*, vol. 119 of *Internat. Ser. Numer. Math.*, pp. 371–381, Birkhauser Boston, Boston, Mass, USA, 1994.

[41] R. S. Maier, "P-symbols, Heun identities, and 3 F2 identities," in *Special Functions and Orthogonal Polynomials*, vol. 471 of *Contemp. Math.*, pp. 139–159, Amer. Math. Soc., Providence, RI, 2008.

[42] K. Takemura, "Heun's equation, generalized hypergeometric function and exceptional Jacobi polynomial," *Journal of Physics A: Mathematical and General*, vol. 45, no. 8, 085211, 14 pages, 2012.

[43] A. V. Shanin and R. V. Craster, "Removing false singular points as a method of solving ordinary differential equations," *European Journal of Applied Mathematics*, vol. 13, no. 6, pp. 617–639, 2002.

[44] A. Ishkhanyan, "The third exactly solvable hypergeometric quantum-mechanical potential," *EPL (Europhysics Letters)*, vol. 115, no. 2, 2016.

[45] A. M. Ishkhanyan and A. E. Grigoryan, "Fifteen classes of solutions of the quantum two-state problem in terms of the confluent Heun function," *Journal of Physics A: Mathematical and General*, vol. 47, no. 46, 465205, 22 pages, 2014.

[46] L. Carlitz, "Some orthogonal polynomials related to elliptic functions," *Duke Mathematical Journal*, vol. 27, pp. 443–459, 1960.

[47] K. Kuiken, "Heun's equation and the hypergeometric equation," *SIAM Journal on Mathematical Analysis*, vol. 10, no. 3, pp. 655–657, 1979.

[48] J. N. Ginocchio, "A class of exactly solvable potentials. I. One-dimensional Schrödinger equation," *Annals of Physics*, vol. 152, no. 1, pp. 203–219, 1984.

[49] G. Valent, "An integral transform involving Heun functions and a related eigenvalue problem," *SIAM Journal on Mathematical Analysis*, vol. 17, no. 3, pp. 688–703, 1986.

[50] G. S. Joyce, "On the cubic lattice Green functions," *Proceedings of the Royal Society A Mathematical, Physical and Engineering Sciences*, vol. 445, no. 1924, pp. 463–477, 1994.

[51] G. S. Joyce and R. T. Delves, "Exact product forms for the simple cubic lattice Green function. I," *Journal of Physics A: Mathematical and General*, vol. 37, no. 11, pp. 3645–3671, 2004.

[52] G. S. Joyce and R. T. Delves, "Exact product forms for the simple cubic lattice Green function. II," *Journal of Physics A: Mathematical and General*, vol. 37, no. 20, pp. 5417–5447, 2004.

[53] R. S. Maier, "On reducing the Heun equation to the hypergeometric equation," *Journal of Differential Equations*, vol. 213, no. 1, pp. 171–203, 2005.

[54] M. van Hoeij and R. Vidunas, "Belyi functions for hyperbolic hypergeometric-to-Heun transformations," *Journal of Algebra*, vol. 441, pp. 609–659, 2015.

[55] R. Vidunas and G. Filipuk, "Parametric transformations between the Heun and Gauss hypergeometric functions," *Funkcialaj Ekvacioj. Serio Internacia*, vol. 56, no. 2, pp. 271–321, 2013.

[56] R. Vidunas and G. Filipuk, "A classification of coverings yielding Heun-to-hypergeometric reductions," *Osaka Journal of Mathematics*, vol. 51, no. 4, pp. 867–903, 2014.